\documentclass[a4paper,12pt]{article}

\begin{document}
\author{Sergey V. Ludkovsky.}
\title{Microbundles over topological rings.}
\date{29 January 2019}
\maketitle
\begin{abstract}
The article is devoted to microbundles over topological rings. Their
structure, homomorphisms, automorphisms and extensions are studied.
Moreover, compactifications and inverse spectra of microbundles over
topological rings are investigated. \footnote{ 
\par keywords: microbundle; topological ring; homomorphism; inverse spectrum}
\end{abstract}

\par address: Dep. Appl. Mathematics, Moscow State Techn. University MIREA,
av. Vernadsky 78, Moscow, 119454, Russia \par  e-mail:
sludkowski@mail.ru

\section{Introduction.}
Microbundles being generalizations of topological manifolds and
geometric bundles compose a large area in topology and algebraic
topology \cite{fedchigb, milnor64,sulwintb}. Though microbundles on
${\bf R}^n$ and Banach spaces over $\bf R$ were investigated, but on
topological modules over rings they were not broadly studied. On the
other hand, topological rings and topological fields other than $\bf
R$ and $\bf C$ are important not only in algebraic topology and
general topology, but also in their applications (see, for example,
\cite{archtkachb,bourbalg,ludqimb}-\cite{luddfnafjms9,roo,sch1,weil,wieslb}
and references therein). Investigations of microbundles over
topological rings are motivated by problems of general topology,
algebraic topology, algebraic geometry, representation theory,
bundles over topological groups and group rings \cite{felldorb},
mathematical physics.
\par This article is devoted to investigations
of microbundles over topological rings. Relations between
microbundle and manifold structures are elucidated in Theorems 19,
20 and Corollaries 21, 23-25. In Lemmas 14, 15 and Proposition 16
extensions of microbundles are studied. Inverse spectra of
microbundles are investigated in Theorems 27, 29, 30, Propositions
31 and 38. Extensions and homomorphisms of microbundles related with
extensions and homomorphisms of topological modules and topological
rings are described in Theorem 32 and Corollary 33.
Compactifications of microbundles are studied in Theorem 34 and
Corollary 35.
\par All main results of this paper are obtained for the first time.
They can be used for further studies not only in topology and its
applications such as geometry, algebraic topology, representation
theory and mathematical analysis, but also in mathematical physics.

\section{Microbundles.}
\par {\bf 1. Notation.} Let $\bf F$ be a topological ring
such that its topology $\tau _{\bf F}$ is neither discrete nor
antidiscrete. Then $({\bf F}, \tau _{\bf F})$ is called a proper
topological ring. We consider a topological left module $X_{\bf F}$
over $\bf F$, or shortly $X$ if $\bf F$ is specified. The ring $\bf
F$ is supposed to be associative and commutative relative to the
addition, but may be noncommutative or nonassociative relative to
the multiplication.
\par Henceforward, it will also be written shortly a ring or a module
instead of a topological ring or a topological module. Their
homomorphisms will be supposed being continuous. Neighborhoods in
topological spaces, modules, rings will be open if something other
will not be specified and the topological terminology is used in the
sense of the book \cite{eng}.

\par {\bf 2. Definition.} Suppose the following. \par $(2.1)$ There are are topological
spaces $A$ and $E$.
\par $(2.2)$ There are continuous maps $i: A\to E$ and $p: E\to A$ such
that $p\circ i=id$, where $id: A\to A$ is the identity map, $p\circ
i$ denotes a composition of maps. Then $A$ will be called a base
space, $E$ will be called a total space.
\par $(2.3)$ For each $b\in A$ neighborhoods $U$ of $b$ and $V$ of $i(b)$ exist
such that $i(U)\subset V$ and $p(V)\subset U$ and $V$ is
homeomorphic to $U\times X$, where $h_V: V\to U\times X$ is a
homeomorphism. \par $(2.4)$ There are continuous maps: a projection
$\hat{\pi }_1: U\times X\to U$, an injection $\iota _0: U\to U\times
X$, also a projection $\hat{\pi }_2: U\times X\to X$ such that
$\hat{\pi }_1(d,x)=d$, $ ~ \iota _0(d)=(d,0)$ and $\hat{\pi
}_2(d,x)=x$ for each $d\in U$ and $x\in X$. They are supposed to
satisfy the identity: $\hat{\pi }_1\circ \iota _0|_U=p|_V \circ i|_U
$, where $i|_U$ denotes the restriction of $i$ to $U$.
\par If the conditions $(2.1)-(2.4)$ are satisfied, then it will be said that
they define a microbundle ${\cal B}={\cal B}(A, E, {\bf F}, X, i,
p)$ with a fibre $X=X_{\bf F}$ of ${\cal B}$. \par If a fibre $X$ is
finite dimensional over the ring $\bf F$, that is $X={\bf F}^n$ with
$n\in \bf N$, then $n$ is called the fibre dimension of ${\cal B}$
over ${\bf F}$. If $X$ is infinite dimensional over ${\bf F}$, then
it is said that the microbundle $\cal B$ has an infinite fibre
dimension over ${\bf F}$. If some data are specified, like $\bf F$
or $X$, they can be omitted from ${\cal B}(A, E, {\bf F}, X, i, p)$
in order to shorten the notation.

\par {\bf Examples. 3.1.} In particular, if $E=A\times X~$, $i=\iota _0~$, $p=\hat{\pi }_1$,
then such a microbundle will be called the standard trivial
microbundle and it will be denoted by ${\sf s}_{A,X}$ or ${\sf s}$.
\par {\bf 3.2.} Suppose that $\xi $ is a vector bundle over $A$ with
a fibre $X$ over a field $\bf F$ and a structure group $GL(X)$ of
all continuous linear automorphisms $T: X\to X$. Suppose also that
$E$ is its total space, $p: E\to A$ is a projection, $i: A\to E$ is
a zero cross-section. This provides an underlying microbundle $|\xi
|$ of $\xi $.

\par {\bf 4. Definitions.} Let ${\cal B}_1={\cal B}(A_1,E_1,X_1,i_1,p_1)$
and ${\cal B}_2={\cal B}(A_2,E_2,X_2,i_2,p_2)$ be two microbundles.
Let also neighborhoods $V_1$ of $i_1(A_1)$ in $E_1$ and $V_2$ of
$i_2(A_2)$ in $E_2$ and homeomorphisms $g: V_1\to V_2$ and $s:
A_1\to A_2$ exist such that $p_2\circ g|_{V_1}=s\circ {p_1}|_{V_1}$
and $g\circ i_1=i_2\circ s$ and $p_2\circ i_2\circ s=s\circ p_1\circ
i_1$. Then these microbundles are called base neighbor isomorphic.
For short it will be said "isomorphic microbundles" instead of "base
neighbor isomorphic microbundles".
\par A microbundle is called trivial if it is isomorphic to the
standard trivial microbundle ${\sf s}$. We remind also the
following. \par Let $M$ be a topological space and let $X$ be a left
module over a proper topological ring $\bf F$. Suppose also that
\par $(4.1)$ $M$ has a covering $ \{ U_k: ~ k\in K \} $, that is
$\bigcup_{k\in K} U_k=M$, where $U_k$ is an open subset in $M$ for
each $k\in K$, where $K$ is a set; \par $(4.2)$ for each $k\in K$
there exists a homeomorphism $\phi _k: U_k\to V_k$, where $V_k$ is
open in $X$.
\par Then $M$ is called a topological manifold on $X$. A triple
$(U_k, \phi _k, V_k)$ is called a chart, where $k\in K$. A
collection $\{ (U_k, \phi _k, V_k): k\in K \} $ of charts is called
an atlas of $M$ and denoted by $At (M)$.

\par {\bf 5. Lemma.} {\it Suppose that $M$ is a topological manifold
on $X$ possessing an atlas of $M$ with charts homeomorphic to $X$
and supplied with a diagonal map $\Delta : M\to M\times M$. Then
this manifold induces a microbundle with $A=M$, $E=M\times M$,
$i=\Delta $.}
\par {\bf Proof.} Evidently $\hat{\pi }_1\circ \Delta =id$ on $M$.
For any $m\in M$ take a neighborhood $U$ such that a homeomorphism
$f: U\to X$ exists. It induces a map $g: U\times U\to U\times X$
such that $g(x,y)=(x,f(y)-f(x))$. Then $g$ is a homeomorphism of
$U\times U$ onto $U\times X$ such that $\hat{\pi }_1\circ \Delta
|_U= \hat{\pi }_1\circ g$ and $g\circ \Delta |_U={\iota _0}|_U$ and
${\hat{\pi }_1}|_{U\times U}=\hat{\pi }_1\circ g$.
\par {\bf 6. Definition.} The microbundle provided by Lemma 5 is called
the tangent microbundle of $M$ and denoted by ${\sf t}_M$ or ${\sf
t}$.
\par {\bf 7. Note.} We consider a local field $\bf K$.
This is a finite algebraic extension of the field $\bf Q_p$ of
$p$-adic numbers with a multiplicative nontrivial norm $|\cdot
|_{\bf K}$ extending that of $\bf Q_p$ \cite{weil,wieslb}.
\par Using antiderivation operators in the sense of Schikhof
\cite{sch1} in Section 2 in \cite{ludspnamijmms4} and in Section 3
in \cite{ludladijmms5} were defined and investigated manifolds over
$\bf K$ of classes $\mbox{ }_PC_0((t,s))$ and $\mbox{
}_S^lC^{(q+l,n-1)}$ respectively.
\par It appears that for them tangent microbundle structures also exist.

\par {\bf 8. Theorem.} {\it Let $M$ be over ${\bf F}=\bf K$ either a
$\mbox{ }_PC_0((t,s))$-manifold with $s\ge 2$ and $t\ge 0$ or a
$\mbox{ }_S^lC^{(q+l,n-1)}$-manifold  with $l\ge 2$, $q\ge 0$ and
$n\ge 1$ (see Note 7). Let also $TM$ be its tangent vector bundle.
Then the underlying microbundle $|TM|$ corresponding to $TM$ is
isomorphic to the tangent microbundle ${\sf t}_M$ of $M$.}
\par {\bf Proof.} In the first case in view of Theorem 2.7 in \cite{ludspnamijmms4} a
clopen neighborhood ${\tilde T}M$ of $M$ in $TM$ exists together
with an exponential $C_0((t,s))$-mapping $\exp : {\tilde T}M\to M$
of ${\tilde T}M$ on $M$. In the second case by virtue of Theorem
3.23 in \cite{ludladijmms5} there exist a clopen neighborhood
${\tilde T}M$ of $M$ in $TM$ and an exponential $\mbox{
}_S^lC^{(q+l,n-1)}$-mapping $\exp : {\tilde T}M\to M$ of ${\tilde
T}M$ on $M$. This mapping $\exp $ is induced by considering
geodesics in $M$ over $\bf K$. \par There exists the natural
embedding $\phi: M\to {\tilde T}M$ such that $M\ni y\mapsto \phi
(y)=(x,0)\in {\tilde T}M$. Therefore, a map $f: {\tilde T}M\to
M\times M$ can be defined for which \par $f(x,v)=(x, \exp _x(v))$
for each $(x,0)\in \phi (M)$ and $v\in {\tilde T}_xM$. We apply to
$f$ a non-archimedean analog of the Taylor Theorem A.1 in
\cite{luddfnafjms9} and the inverse function Theorem A.4 in
\cite{ludqimb,ludqimnafjms5} (see also \cite{sch1}). They imply that
for each $(x,0)\in \phi (M)$ a neighborhood $U_x$ of $(x,0)$ in
${\tilde T}M$ exists such that the restriction $f|_{U_x}$ is a
diffeomorphism on a neighborhood $V_y$ of $(y,y)\in M\times M$,
where $y\in M$ and $\phi (y)=(x,0)$. Taking a covering of the
diagonal $D_M=\{ (y,y): y\in M \} $ in $M\times M$ by such
neighborhoods $V_y$ provides neighborhoods $V$ of $D_M$ in $M\times
M$ and $U$ of $\phi (M)$ in ${\tilde T}M$ such that $f: U\to V$ is a
diffeomorphism of $U$ onto $V$. Therefore, from the commutative
diagram of this situation it follows that the underlying microbundle
$|TM|$ corresponding to $TM$ is isomorphic to the tangent
microbundle ${\sf t}_M$ of $M$.

\par {\bf 9. Definition.} Let $\bf F$ be a unital topological ring.
\par A topological space $A$ will be called $\bf F$ completely regular, if
it is $T_1$ and for each closed subset $V$ in $A$ and each point
$a\in A\setminus V$ a continuous function $f: A\to \bf F$ exists
such that $f(a)=0$ and $f(V)= \{ 1 \} $.

\par {\bf 10. Realization of trivial microbundles.}
\par Let $\bf F$ be an infinite field with \par $(10.1)$ a topology
induced by a multiplicative norm $|\cdot |_{\bf F}$, where the norm
takes values in $[0, \infty )= \{ t\in {\bf R}: ~ t\ge 0 \} $ and
\par $(10.2)$ let $\Gamma _{\bf F}$ be dense in $(0, \infty )$, \\ where
$\Gamma _{\bf F} = \{ |b|_{\bf F}: 0\ne b \in {\bf F} \} $; let also
\par $(10.3)$ let also $\bf F$ be of zero characteristic $char ({\bf F})=0$.
\par Suppose that $X$ is a Banach space
over the field $\bf F$ with a nontrivial norm taking values in
$\Gamma _{\bf F}\cup \{ 0 \} $. Suppose also that ${\sf s}_{A,X}$ is
a trivial microbundle with a paracompact $\bf F$ completely regular
base space $A$ of zero dimension $dim (A)=0$ and a fibre $X$

\par {\bf 11. Proposition.} {\it If conditions of subsection 10 are
satisfied, then an open subset $U_0$ in ${\sf s}_{A,X}$ exists such
that it is homeomorphic to $A\times X$. Moreover, this homeomorphism
is compatible with the injection and projection maps.}
\par {\bf Proof.} Using the definition of the trivial microbundle
${\sf s}_{A,X}$ we consider an open subset $U$ of $A\times X$.
\par For each closed subset $V$ in $A$ and each point $a\in A\setminus
V$ open neighborhoods $W_a$ of $a$ and $W_V$ of $V$ exist which are
disjoint $W_a\cap W_V=\emptyset $, since a continuous function $f:
A\to \bf F$ exists such that $f(a)=0$ and $f(V)= \{ 1 \} $. Indeed,
one can take $W_a=f^{-1}( \{ b\in {\bf F}: |b|<r_1 \} )$ and
$W_V=f^{-1}( \{ b\in {\bf F}: r_2<|b| \} )$, where $0<r_1<r_2<1$,
$r_1$ and $r_2$ belong to $\Gamma _{\bf F} $. An existence of such
$r_1$ and $r_2$ follows from the condition $(10.2)$, since $|0|_{\bf
F}=0$ and $|1|_{\bf F}=1$ and $|b|_{\bf F}>0$ for each $b\ne 0$.
Thus $A$ is a $T_3$ space.
\par For each $a\in A$ a radius $0<r(a)<\infty $
with $r(a)\in \Gamma _{\bf F}$ exists such that $(a,x)\in U$ for
each $x\in X$ with $|x|_X<r(a)$. Therefore using the base of the
topology in the product $A\times X$ we infer that an open
neighborhood $W_a$ of $a$ exists such that $\rho (a) := \inf \{
r(b): b\in W_a \} >0$, since $U$ is open in $A\times X$. Thus a
covering ${\cal V}=\{ V_a: ~ a\in A \} $ of $U_1$ exists with $V_a=
\{ (b,x)\in U: b\in W_a, x\in X, |x|_X< r(b) \} $, where
$U_1=\bigcup_{a\in A}V_a$ is a proper open subset in $A\times X$.
Since $A$ is paracompact, this covering ${\cal V}$ of $U_1$ contains
a subcovering ${\cal W}\subset \cal V$ such that ${\cal W} = \{ V_a:
~ a\in A_0 \} $ with $A_0\subset A$ and $ \{ W_a: a\in A_0 \} $ is a
locally finite covering of $A$.
\par The topological space $A$ is normal, since it is $T_1\cap T_3$ and
paracompact (see Section 1.5 and Theorem 5.1.5 in \cite{eng}). Let
${\cal P}=\{ P_j: ~ j\in J \} $ be an open locally finite covering
of $A$, where $J$ is some set, $A=\bigcup_{j\in J} P_j$. From the
lemma about shrinking of covering (see Lemma 5.1.6 in \cite{eng}) it
follows that it contains a covering ${\cal C}=\{ C_j: ~ j\in J \} $
by closed subsets $C_j$ such that $C_j\subset P_j$ for each $j\in
J$; $ ~ A=\bigcup_{j\in J} C_j$. On the other hand, the topological
space $A$ is zero-dimensional, consequently, each subset $C_j$ can
be chosen clopen (closed and open simultaneously) in $A$ (see
Sections 6.2 and 7.1 in \cite{eng}).
\par For each $C_j$ a continuous function $f_j: X\to {\bf F}$ exists
such that $f_j(a)=0$ for each $a\in A\setminus C_j$ and $f_j(a)=1$
for each $a\in C_j$. Then we take a function $f(a)=\sum_{j\in
J}f_j(a)$. Since $char ({\bf F})=0$ and the covering $\cal C$ is
locally finite, then $0<|f(a)|<\infty $ for each $a\in A$. This
implies that a function $g_j(a)=f_j(a)/f(a)$ is continuous for each
$j\in J$ and their sum $g(a)=\sum_{j\in J}g_j(a)=1$ is unit for each
$a\in A$. Thus a family $\{ g_j: j\in J \} $ is the partition of
unity for the covering $\cal P$.
\par Then we consider balls $B({\bf F},t_0,r)=  \{ t\in {\bf F}: |t-t_0|_{\bf F}\le r
\} $ in ${\bf F}$, where $t_0\in {\bf F}$, $ ~ 0<r<\infty $.
Applying this partition of unity to $\cal W$ we get a continuous
function $h: A\to (B({\bf F},0,1)\setminus \{ 0 \})$ such that if
$(a,x)\in A\times X$ and $|x|_X<|h(a)|_{\bf F}$, then $(a,x)\in U$,
where $~|h(a)|_{\bf F}\ge \min (\rho (a),1) r_2$, where $r_2\in
\Gamma _{\bf F}$ with $0<r_2<1$. For each $x\in X$ it is possible to
choose $\xi (x)\in \bf F$ such that $|x|_X=|\xi (x)|_{\bf F}$, since
$|x|_X\in \Gamma _{\bf F}\cup \{ 0 \} $. We put $\psi (a,x)=(a,
(h(a)-\xi (x))^{-1}x)$ for each $(a,x) \in U_0$, where $U_0= \{
(a,x): a\in A, x\in X, |x|_X<|h(a)|_{\bf F} \} $.
\par Since $\Gamma _{\bf F}$ is dense in $(0,\infty )$, then for
each $\epsilon >0$ and $a\in A$ a vector $x\in X$ exists such that
$(a,x)\in U_0$ and $|h(a)-\xi (x)|_{\bf F}<\epsilon $. Therefore
$\psi : U_0\to A\times X$ is a homeomorphism of an open proper
subset $U_0$ of ${\sf s}_{A,X}$ onto ${\sf s}_{A,X}$.

\par {\bf 12. Definition.} For a microbundle
${\cal B}(A, E, {\bf F}, X, i, p)$ and a topological space $A_1$ and
a continuous map $f: A_1\to A$ an induced microbundle ${\cal B}(A_1,
E_1, {\bf F}, X, i_1, p_1)$ is defined with a total space $E_1 = \{
(a_1,e)\in A_1\times E: f(a_1)=p(e) \} $, where
$i_1(a_1)=(a_1,i(f(a_1)))$ for each $a_1\in A_1$, $~p_1(a_1,e)=a_1$
for each $(a_1,e)\in E_1$. The induced microbundle ${\cal B}(A_1,
E_1, {\bf F}, X, i_1, p_1)$ is also denoted by $f^*{\cal B}(A, E, X,
i, p)$.
\par Particularly if $f$ is an inclusion map of $A_1$ into $A$, then
$f^*{\cal B}(A, E, {\bf F}, X, i, p)$ is a so called restricted
microbundle ${\cal B}(A, E, {\bf F}, X, i, p)|_{A_1}={\cal
B}(A_1,E_2,{\bf F},X,i_2,p_2)$ with $E_2=p^{-1}(A_1)$,
$~i_2=i|_{A_1}$, $~p_2=p|_{E_2}$.

\par  {\bf 13. Cone over $\bf F$.}
\par Assume that \par $(13.1)$ $\bf F$ is an infinite unital ring with
a topology induced by a nontrivial norm taking values in $[0,\infty
)$.
\par We put $CA_1=(A_1\times B({\bf F},0,1))/ (A_1\times \{ 0 \} )$
to be a cone of a topological space $A_1$ over the ring $\bf F$,
where $B({\bf F},x,r) := \{ y\in {\bf F}: ~ |x-y|_{\bf F}\le r \}$,
$x\in \bf F$, $0<r<\infty $.
\par $(13.2)$ For topological spaces $A$ and $A_1$ and a continuous map
$f: A_1\to A$ \\ let $A\bigcup_f CA_1 = (A\cup CA_1)/\Xi $ be a
mapping cone of $f$, where $\Xi =\Xi _f$ denotes an identification
$(a_1,1)\Xi f(a_1)$ for each $a_1\in A_1$.

\par {\bf 14. Lemma.} {\it If Conditions $(13.1)$ and $(13.2)$ are
fulfilled and a ring $\bf F$ is path-connected and a microbundle
${\cal B}(A, E, {\bf F}, X, i, p)$ can be extended to a microbundle
over $A\bigcup_f CA_1$, then $f^*{\cal B}(A, E, {\bf F}, X, i, p)$
is trivial.}
\par {\bf Proof.} Apparently the composition
${A_1}_{\overrightarrow{f}} A_{\overrightarrow{q}} A\bigcup_f CA_1$
is null-homotopic, since $\bf F$ is path-connected, where $q$ is an
embedding of $A$ into $A\bigcup_f CA_1$. Consequently, $f^*{\cal
B}(A, E, {\bf F}, X, i, p)$ is trivial.

\par {\bf 15. Lemma.} {\it If Conditions $(10.1)-(10.3)$ and $(13.2)$
are satisfied and an induced microbundle $f^*{\cal B}(A, E, {\bf F},
X, i, p)$ is trivial, then ${\cal B}(A, E, {\bf F}, X, i, p)$ can be
extended over $A\bigcup_f CA_1$.}
\par {\bf Proof.} At first we take the mapping cylinder $Z=A\bigcup_f(A_1\times
B({\bf F},0,1))$  of $f$, where $A\bigcup_f(A_1\times B({\bf
F},0,1))= [A\bigcup(A_1\times B({\bf F},0,1))]/\Xi _f$. Then the
microbundle ${\cal B}(A, E, {\bf F},  X, i, p)$ can be extended to a
microbundle ${\cal B}_1$ over $Z$, since $A$ is a retract of $Z$.
Therefore the restriction ${\cal B}_1|_{A_1\times \{ 0 \} }$ is
trivial as well, consequently, ${\cal B}_1|_{A_1\times B({\bf
F},0,r)}$ is trivial for each $0<r<1$ with $r\in \Gamma _{\bf F}$.
We fix such $r$.
\par By virtue of Proposition 11 an open subset $U_0$ of the total space
$E_{1,r}$ of the restricted microbundle ${\cal B}_1|_{A_1\times
B({\bf F},0,r)}$ is homeomorphic to $A_1\times B({\bf F},0,r)\times
X$ with a homeomorphism $h$ compatible with injections and
projections. Since $A\bigcup_f CA_1=Z/(A_1\times \{ 0 \} )$, then it
induces from ${\cal B}_1$ a microbundle ${\cal B}_2$ over
$A\bigcup_f CA_1$. It remains to note that a total space $E_2$ of
${\cal B}_2$ is obtained from $E_1$ by an identification
$h^{-1}(A_1\times \{ 0 \} \times x)$ with $x$ for each $x\in X$,
where $E_1$ is a total space of ${\cal B}_1$.

\par {\bf 16. Proposition.} {\it If Conditions $(13.1)$ and $(13.2)$ are
satisfied and a ring $\bf F$ is zero-dimensional, $dim ({\bf F})=0$,
then a microbundle ${\cal B}(A, E, {\bf F}, X, i, p)$ can be
extended to a microbundle over $A\bigcup_f CA_1$.}
\par {\bf Proof.} At first we take a partition of the unit ball
$B({\bf F},0,1)$ into two disjoint clopen subsets $K_0$ and $K_1$
such that $0\in K_0$ and $1\in K_1$, that is $K_0\cap K_1=\emptyset
$ and $K_0\cup K_1=B({\bf F},0,1)$, since $dim ({\bf F})=0$.
Therefore $A\bigcup_f CA_1$ is the disjoint union of two clopen
subsets $A_2:=[A\bigcup (A_1\times K_1)]/\Xi _f$ and $A_3:=[A\bigcup
(A_1\times K_0)]/(A_1\times \{ 0 \} )$. \par Let $T: X\to X$ be any
left $\bf F$ linear automorphism of a topological left module $X$
over $\bf F$ (that is $T$ and $T^{-1}$ are continuous), $V=V_b$ be a
neighborhood of $b=i\circ f(a_1)$ in $E$, $h_V: V\to U\times X$ be a
homeomorphism, $U$ be a neighborhood of $f(a_1)$ in $A$, $\hat{\pi
}_2: U\times X\to X$ be a projection such that $\hat{\pi }_2(d,x)=x$
for each $d\in U$ and $x\in X$ (see Definition 2). We put $E_2=[E
\bigcup (A_1\times K_1\times X)]/\Xi _{g}$, where $g: A_1\times X\to
E$ is a continuous mapping such that $p\circ g(a_1,x)=i\circ f(a_1)$
and $\hat{\pi }_2\circ h_V\circ g(a_1,x)=Tx$ for each $a_1\in A_1$
and $x\in X$, where $\Xi _g$ identifies $(a_1,1,x)$ with $g(a_1,x)$.
\par Take any automorphism $g_2: B({\bf F},0,1)\to B({\bf F},0,1)$
such that $g_2(1)=1$ and $g_2(0)=0$ (that is $g_2$ and $g^{-1}_2$
are continuous). An injection $i: A\to E$ has a continuous extension
$i: A_2\to E_2$ such that $p\circ i(a_1,t)=(i\circ f(a_1),g_2(t))$
and $\hat{\pi }_2\circ h_V\circ i(a_1,t)=(t+\beta (1-t)) \hat{\pi
}_2\circ h_V\circ i\circ f(a_1)$ for each $t\in K_1\setminus \{ 1 \}
$ and $a_1\in A_1$, where $V$ is a neighborhood of $i(a_1,t)$ in
$E_2$, $~\beta $ is a fixed element in $\bf F$. Therefore the
projection $p: E\to A$ has a continuous extension on $E_2$ such that
$p: E_2\to A_2$ with $p\circ h_V^{-1}(b,t,x)= (a_1,g_2^{-1}(t))$ for
each $t\in K_1\setminus \{ 1 \} $ and $a_1\in A_1$ with $b=i\circ
f(a_1)$. That is $p\circ i=id$ on $A_2$. Thus a microbundle ${\cal
B}_2={\cal B}(A_2,E_2, {\bf F}, X, i, p)$ is an extension of ${\cal
B}(A, E, {\bf F}, X, i, p)$.
\par On $A_3$ a microbundle ${\cal B}_3={\cal B}(A_3, E_3, {\bf F}, X, i_3, p_3)$ exists,
which may be in particular trivial ${\sf s}_{A_3,X}$.
\par For mappings $f_j: B_j\to C_j$ for each $j\in \{ 1, 2 \} $
and $B=B_1\cup B_2$ and $C=C_1\cup C_2$ with $B_1\cap B_2=\emptyset
$ and $C_1\cap C_2=\emptyset $ by $f_1\nabla f_2$ is denoted their
combination such that $(f_1\nabla f_2)(b_j)=f_j(b_j)$ for each
$b_j\in B_j$ and $j\in \{ 1,2 \} $. Therefore the combination ${\cal
B}_2\nabla {\cal B}_3= {\cal B}(A_2\cup A_3,E_2\cup E_3,{\bf F}, X,
i_2\nabla i_3, p_2\nabla p_3)$ of microbundles ${\cal B}_2$ and
${\cal B}_3$ provides the extension over $A\bigcup_f CA_1$ of the
microbundle ${\cal B}(A, E, {\bf F}, X, i, p)$.

\par {\bf 17. Definition.} Suppose that $M_1$ and
$M$ are topological manifolds (see Definition 4) on ${\bf F}^{m_1}$
and ${\bf F}^m$ over a topological ring $\bf F$ such that
$M_1\subset M$, where $m_1$ and $m$ are cardinals such that $m_1\le
m$, where ${\bf F}^m$ is supplied with the Tychonoff product
topology, ${\bf F}$ has the topological weight $\tau = w{\bf F}\ge
\aleph _0$. Suppose also that a neighborhood $U$ of $M_1$ in $M$
exists and a retraction $p: U\to M_1$ is such that
${M_1}_{\overrightarrow{i}}U_{\overrightarrow{p}}M_1$ forms a
microbundle, where $i$ is an inclusion map. Then it will be said
that $M_1$ has a microbundle neighborhood in $M$. It will be denoted
by ${\cal N}={\cal N}(M_1,M,i,p)$.
\par In particular, if $U$ and $p$ can be chosen such that
${\cal N}$ is trivial, then $M_1$ has a product neighborhood.
\par {\bf 18. Corollary.} {\it Assume that Conditions $(10.1)-(10.3)$
are satisfied. Then in the notation of Definition 17 $M_1$ has a
trivial microbundle neighborhood if and only if a neighborhood $U$
of $M_1$ in $M$ exists such that the pair $(U,M_1)$ is homeomorphic
to $M_1\times (X,0)$.}
\par This follows from Proposition 11.
\par {\bf 19. Theorem.} {\it If $M_1$ and $M$ are topological
manifolds over a topological ring $\bf F$ of zero small inductive
dimension, and $M_1$ is a closed submanifold in $M$, then a
retraction $\check{r}: M\to M_1$ exists.}
\par {\bf Proof.} Since $ind ({\bf F})=0$, then the small inductive
dimension of ${\bf F}^m$ is zero, $ind ({\bf F}^m)=0$ (see Theorem
6.2 and Section 7.1 in \cite{eng}). An atlas of $M$ has a refinement
being an atlas with charts homeomorphic to clopen subsets in ${\bf
F}^m$, since ${\bf F}^m$ has a base of its topology consisting of
clopen subsets. Thus we consider that each chart $V_j$ of the atlas
$At (M) = \{ (V_j, f_j): j\in J \} $ of $M$ is homeomorphic to a
clopen subset $W_j=f_j(V_j)$ in ${\bf F}^m$, where $J$ is a set. A
similar choice of an atlas $At (M_1)$ can be made for $M_1$. \par
Note that ${\bf F}^m$ has an embedding $g$ into the generalized
Cantor discontinuum $D^{\tau m}$, because $\tau \ge \aleph _0$.
Therefore, $M$ has an embedding into $D^t$, where $t= \tau mn$,
$\tau $ is a topological weight of $\bf F$, $n$ is a cardinality of
$J$. We take the closure $C=cl (g(M))$ of $g(M)$ in $D^t$. There
exists a retraction $q: C\to C_1$, where $C_1=cl (g(M_1))$ (see
\cite{eng,fedchigb,fedumn86}). Therefore, the restriction
$q|_{g(M)}$ induces $\check{r}=g^{-1}\circ q\circ g|_V$, where
$V=\check{r}^{-1}(M_1)=M$, since $M_1$ is closed in $M$; $ ~g(M_1)$
is dense in the compact space $C_1$ and $q(c_1)=c_1$ for each
$c_1\in C_1$.
\par {\bf 20. Theorem.} {\it Assume
that manifolds $M$ and $M_1$ are both over a field ${\bf F}=\bf K$
either $\mbox{ }_PC_0((t,s))$-manifolds with $s\ge 2$ and $t\ge 0$
or $\mbox{ }_S^lC^{(q+l,n-1)}$-manifolds with $l\ge 2$, $q\ge 0$ and
$n\ge 1$ (see Definitions 6, 17 and Note 7), where $M_1$ is a closed
submanifold in $M$. Then the total space of $i^*{\sf t}_M$ is
homeomorphic to the total space of $\check{r}^*{\sf t}_{M_1}$.}
\par {\bf Proof.} By virtue of Theorem 19 we consider the case
when there is a retraction $\check{r}: M\to M_1$. The total space
$E=E(i^*{\sf t}_M)$ of $i^*{\sf t}_M$ consists of all pairs
$(m_1,(m,k))\in M_1\times M^2$ such that $i(m_1)=m$, consequently,
it is homeomorphic to $M_1\times M$. Then $E_1=E(\check{r}^*{\sf
t}_{M_1})$ is a subspace in $M\times M_1^2$ composed of all pairs
$(m,(m_1,k))$ with $\check{r}(m)=m_1$, hence the total space $E_1$
is homeomorphic to $M\times M_1$.
\par {\bf 21. Corollary.} {\it If the conditions of Theorem 20 are
accomplished, then the submanifold $M_1$ has a microbundle
neighborhood ${\cal N}$ homeomorphic with $i^*{\sf t}_M$.}
\par {\bf 22. Definition.} If a microbundle ${\sf t}_M$
is trivial, then a manifold $M$ is called topologically
parallelizable.
\par {\bf 23. Corollary.} {\it If conditions of Theorem 20
are fulfilled and $M_1$ is topologically parallelizable, then
$M_1\times \{ 0 \} $ has a microbundle neighborhood in $M\times X$
with a normal microbundle ${\cal N}$ homeomorphic to $i^*{\sf
t}_M$.}
\par {\bf Proof.} Since ${\sf t}_{M_1}$ is trivial, then
$\check{r}^*{\sf t}_{M_1}$ is trivial as well. Therefore, the total
space $E(\check{r}^*{\sf t}_{M_1})$ is homeomorphic with $M\times X$
of the canonical trivial microbundle ${\sf s}_{M,X}$, where $X={\bf
F}^m$.
\par {\bf 24. Corollary.} {\it Let the conditions of Theorem 20 be
satisfied and let $M$ and $M_1$ be topologically parallelizable.
Then $M_1\times \{ 0 \} $ has a product neighborhood $V$ in $M\times
X$.}
\par {\bf Proof.}  Since ${\sf t}_M$ is trivial, then ${\cal
N}$ is trivial.
\par {\bf 25. Corollary.} {\it Let $M$ be a compact
topologically parallelizable of dimension $1\le m=dim_{\bf
K}M<\infty $ over a field ${\bf F}=\bf K$ either $\mbox{
}_PC^{(q,n)}$-manifolds or $\mbox{ }_S^lC^{(q,n)}$-manifolds with
$q\ge 1$ and $n\ge 0$ (see Definitions 6 and 17).
\par Then a product neighborhood $V$ of $M\times {\bf K}^{2m+1}$
exists such that it can be embedded into ${\bf K}^{3m+1}$ as a
clopen subset.}
\par {\bf Proof.} By virtue of Theorem 3.21 in \cite{ludladijmms5}
there exists a $\mbox{ }_SC^{(q,n)}$ or $\mbox{
}_PC^{(q,n)}$-embedding $\tau : M\hookrightarrow {\bf K}^{2m+1}$
correspondingly. Since $\bf K$ is the local field, it is
zero-dimensional $dim ({\bf K})=0$ (see \cite{eng,roo,weil}). Each
ball $B({\bf K}^n,x,r) = \{ y\in {\bf K}^n: |x-y|_{\bf K}\le r \} $
is clopen in ${\bf K}^n$, where $x\in {\bf K}^n$, $r\in \Gamma _{\bf
K}$, $n\in {\bf N}$. Therefore a product neighborhood $V$ provided
by Corollary 24 can be chosen homeomorphic to a clopen subset in
${\bf K}^{3m+1}$.

\par {\bf 26. Definition.} Let ${\cal B}_k={\cal B}(A_k,E_k, {\bf F}_k, X_k, i_k,p_k)$
be a family of microbundles over topological rings ${\bf F}_k$ with
$k\in \Lambda $, where $\Lambda $ is a directed set. Let also for
each $k\le n$ in $\Lambda $ a homomorphism $\pi ^n_k: {\cal B}_n\to
{\cal B}_k$ be given satisfying the conditions $(26.1)-(26.5)$:
\par $(26.1)$ $\pi ^n_k = (\pi ^{1,n}_{1,k}; \pi ^{2,n}_{2,k}; \pi
^{3,n}_{3,k}; \pi ^{4,n}_{4,k})$ with
\par $(26.2)$ $\pi ^{1,n}_{1,k}: A_n\to A_k$ and  $\pi ^{2,n}_{2,k}:
E_n\to E_k$, \par $\pi ^{3,n}_{3,k}: {\bf F}_n\to {\bf F}_k$ and
$\pi ^{4,n}_{4,k}: X_n\to X_k$, where
\par $(26.3)$ $\pi ^{2,n}_{2,k}\circ i_n=i_k\circ \pi
^{1,n}_{1,k}$ and $p_k\circ \pi ^{2,n}_{2,k}=\pi ^{1,n}_{1,k}\circ
p_n$, \par $(26.4)$ $h_k\circ \pi ^{2,n}_{2,k}\circ
h_n^{-1}(a_n,u_nx_n+v_ny_n)=$\par $ \pi ^{3,n}_{3,k}(u_n)h_k\circ
\pi ^{2,n}_{2,k}\circ h_n^{-1}(a_n,x_n) + \pi ^{3,n}_{3,k}(v_n)
h_k\circ \pi ^{2,n}_{2,k}\circ h_n^{-1}(a_n,y_n)$,
\par $(26.5)$ $\hat{\pi }_{2,k} \circ h_k\circ \pi ^{2,n}_{2,k}=\pi
^{4,n}_{4,k}\circ \hat{\pi }_{2,n} \circ h_n$ and \par $\pi
^{4,n}_{4,k}(u_nx_n+v_ny_n)=\pi ^{3,n}_{3,k}(u_n) \pi
^{4,n}_{4,k}(x_n)+ \pi ^{3,n}_{3,k}(v_n) \pi ^{4,n}_{4,k}(y_n)$ \\
for every $u_n$ and $v_n$ in ${\bf F}_n$, $x_n$ and $y_n$ in $X_n$,
$a_n\in A_n$, where $h_n: V_n\to U_n\times X_n$ is a local
homeomorphism for an open subset $V_n$ in $E_n$ corresponding to an
open neighborhood $U_n$ of a point $a_n$ in $A_n$; where $\pi ^n_n$
is the identity homomorphism if $n=k$ and $\pi ^k_l\circ \pi
^n_k=\pi ^n_l$ for each $l\le k \le n$ in $\Lambda $.
\par Such a family $\{ {\cal B}_n, \pi ^n_k, \Lambda \} $ will be
called an inverse spectrum of microbundles.
\par Let $A=\lim \{ A_n, \pi ^{1,n}_{1,k},
\Lambda \} $ and $E=\lim \{ E_n, \pi ^{2,n}_{2,k}, \Lambda \} $ be
limits of inverse spectra of topological spaces. Let also ${\bf F}=
\lim \{ {\bf F}_n, \pi ^{3,n}_{3,k}, \Lambda \} $ and $X=\lim \{
X_n, \pi ^{4,n}_{4,k}, \Lambda \} $ be limits of inverse spectra of
topological rings and topological left modules respectively.
\par Let ${\cal B}={\cal B}(A,E,{\bf F},X,i,p)$ be a microbundle
such that for each $n\in \Lambda $ a homomorphism $\pi _n: {\cal
B}\to {\cal B}_n$ exists satisfying analogous to $(26.1)$-$(26.5)$
conditions and the following condition: $\pi ^n_k\circ \pi _n = \pi
_k$ for each $k\le n$ in $\Lambda $. Then it will be said that
${\cal B}$ is a limit of the inverse spectrum $ \{ {\cal B}_n, \pi
^n_k, \Lambda \} $ of microbundles.

\par {\bf 27. Theorem.} {\it If $\{ {\cal B}_n, \pi ^n_k, \Lambda \} $
is an inverse spectrum of microbundles, then its limit exists.}
\par {\bf Proof.} The inverse spectrum of microbundles induces
inverse spectra of topological spaces $\{ A_n, \pi ^{1,n}_{1,k},
\Lambda \} $ and $\{ E_n, \pi ^{2,n}_{2,k}, \Lambda \} $. Therefore
there exist topological spaces $A$ and $E$ being their limits $A=
\lim \{ A_n, \pi ^{1,n}_{1,k}, \Lambda \} $ and $E=\lim \{ E_n, \pi
^{2,n}_{2,k}, \Lambda \} $ (see Section 2.5 in \cite{eng} and
\cite{scepumn1976}). For each $k\in \Lambda $ there are projections
$\pi _{1,k}$ from $A$ onto $A_k$ and $\pi _{2,k}$ from $E$ onto
$E_k$.
\par Then from $(26.3)-(26.5)$ it follows that local
homeomorphisms $h_n: V_n\to U_n\times X_n$ are compatible with
inverse spectra of rings $\{ {\bf F}_n, \pi ^{3,n}_{3,k}, \Lambda \}
$ and of left modules $\{ X_n, \pi ^{4,n}_{4,k}, \Lambda \} $, where
$\pi ^{3,n}_{3,k}: {\bf F}_n\to {\bf F}_k$ is a homomorphism of
topological rings and $\pi ^{4,n}_{4,k}: X_n\to X_k$ is a
homomorphism of left topological modules such that
\par $(27.1)$ $\pi ^{4,n}_{4,k}(u_nx_n+v_ny_n)= \pi
^{3,n}_{3,k}(u_n)\pi ^{4,n}_{4,k}(x_n) + \pi ^{3,n}_{3,k}(v_n) \pi
^{4,n}_{4,k}(y_n)$ \\ for every $u_n$ and $v_n$ in ${\bf F}_n$,
$x_n$ and $y_n$ in $X_n$, $~k\le n$ in $\Lambda $. For each $n\in
\bf N$ the left module $X_n$ has also a structure of a commutative
group relative to the addition on it.
\par Therefore there exist a topological ring ${\bf F}=\lim \{ {\bf
F}_n, \pi ^{3,n}_{3,k}, \Lambda \} $ and a commutative group
(relative to the addition) $X=\lim \{ X_n, \pi ^{4,n}_{4,k}, \Lambda
\} $ (see \cite{bourbalg,eng}). There are projections
(homomorphisms) $\pi _{3,k}$ from ${\bf F}$ onto ${\bf F}_k$ and
$\pi _{4,k}$ from $X$ onto $X_k$. Each $u$ in $\bf F$ has the form
$u=(u_n: n\in \Lambda )$ such that $(\forall n\in \Lambda , u_n\in
{\bf F}_n, ~ \forall k\in \Lambda , \forall n\in \Lambda , ~ [(k\le
n)\Rightarrow (\pi ^{3,n}_{3,k}(u_n)=u_k)])$. Each $x\in X$ is of
the form $x=(x_n: n\in \Lambda )$ such that $(\forall n\in \Lambda ,
x_n\in X_n, ~ \forall k\in \Lambda , \forall n\in \Lambda , ~ [(k\le
n)\Rightarrow (\pi ^{4,n}_{4,k}(x_n)=x_k)])$. A base of a topology
on $\bf F$ consists of all subsets $\pi _{3,k}^{-1}(S_k)$ with $S_k$
open in ${\bf F}_k$ and $k\in \Lambda $. Similarly a base of a
topology on $X$ is composed of all subsets $\pi _{4,k}^{-1}(Y_k)$
with $Y_k$ open in $X_k$ and $k\in \Lambda $. Therefore ${\bf F}$
acts continuously on $X$ as $ux=(u_nx_n: n\in \Lambda )$ for each
$u\in \bf F$ and $x\in X$, hence $X$ is a topological left module
over $\bf F$.
\par We have that
$\pi ^{1,n}_{1,k}\circ \pi _{1,n}(a)=\pi _{1,k}(a)$ for each $a\in
A$ and $\pi ^{2,n}_{2,k}\circ \pi _{2,n}(e)=\pi _{2,k}(e)$ for each
$e\in E$ and every $k\le n$ in $\Lambda $. Therefore from Conditions
$(26.1)$ and $(26.2)$ we infer that there exists an injection $i:
A\to E$ such that $i(a)=(b_n: n\in \Lambda )$ satisfying $(\forall
n\in \Lambda , b_n\in E_n, ~ b_n=i_n(a_n) )$ for each $a\in A$,
since $a=(a_n: n\in \Lambda )$ such that $(\forall n\in \Lambda , ~
a_n\in A_n, ~ \forall k\in \Lambda , ~ \forall n\in \Lambda ,
~[(k\le n) \Rightarrow (\pi ^{1,n}_{1,k}(a_n)=a_k )])$. Moreover,
there exists a projection $p: E\to A$ such that $p(b)=(a_n: n\in
\Lambda )$ satisfying the following condition $(\forall n\in \Lambda
, a_n\in E_n, ~ a_n=p_n(b_n) )$ for each $b\in E$, since $b=(b_n:
n\in \Lambda )$ such that $(\forall n\in \Lambda , ~ b_n\in E_n, ~
\forall k\in \Lambda , ~ \forall n\in \Lambda , ~[(k\le n)
\Rightarrow (\pi ^{2,n}_{2,k}(b_n)=b_k )])$. Since $p_n\circ
i_n=id_n$ for each $n\in \Lambda $, then $p\circ i=id$.
\par On the other hand, if $h_k: V_k\to U_k\times X_k$ is a
homeomorphism, where $V_k$ is an open subset in $E_k$ and $U_k$ is
an open subset in $A_k$, then $\pi _{2,k}^{-1}(V_k)=V$ is open in
$E$ and $\pi _{1,k}^{-1}(U_k)=U$ is open in $A$; $ ~ \pi
_{4,k}(X)=X_k$; $~ \pi _{3,k}({\bf F})={\bf F}_k$. Thus bases of
topologies in $A$ and $E$ induce a local homeomorphism $h: V\to
U\times X$ for the corresponding open subsets $V$ in $E$ and $U$ in
$A$, where $h(v)=(a,x)$ with $(a,x)=((a_n,x_n): n\in \Lambda )$ such
that $(\forall n\in \Lambda , (a_n,x_n)=h_n(v_n), a_n\in A_n, x_n\in
X_n)$ for each $v\in V$, where $v=(v_n: n\in \Lambda )$ such that
$(\forall n\in \Lambda , v_n\in E_n, ~ v_n=\pi _{2,n}(v))$. Hence
$\hat{\pi }_1\circ h(v)=p(v)$ for each $v\in V$, since Identity
$(27.1)$ is satisfied and $\hat{\pi }_{1,k}\circ
h_k|_{V_k}=p_k|_{V_k}$ for each $k\in \Lambda $ (see Definition 2).
There is the natural injection $\iota _0: U\hookrightarrow U\times
X$. Then we deduce that $\hat{\pi }_1\circ \iota _0(a)=p|_V \circ
i(a)$ for each $a\in U$, since $\hat{\pi }_{1,k}\circ \iota
_{0,k}|_{U_k}=p_k|_{V_k} \circ i_k|_{U_k}$ for each $k\in \Lambda $.
Thus ${\cal B}(A,E,{\bf F},X,i,p)=\lim \{ {\cal B}_n, \pi ^n_k,
\Lambda \} $.

\par {\bf 28. Definition.} Let ${\bf S}_1=\{ {\cal B}_n, \pi ^n_k, \Lambda \} $
and ${\bf S}_2=\{ {\cal C}_n, \breve{\pi }^n_k, \Upsilon \} $ be two
inverse spectra of micronbundles, let also \par $(28.1)$ $q:
\Upsilon \to \Lambda $ be a map and \par $(28.2)$ $T= \{ \forall
k\in \Upsilon ~ t_k: {\cal B}_{q(k)}\to {\cal C}_k \} $ \\ be a
family of homomorphisms satisfying analogous to $(26.1)$-$(26.5)$
conditions such that for each $k\le n$ in $\Upsilon $ there exists
$m\in \Lambda $ with $m\ge q(n)$ and $m\ge q(k)$ for which the
following identity is satisfied:
\par $(28.3)$ $t_k\circ \pi ^m_{q(k)}= \breve{\pi }^n_k\circ
t_n\circ \pi ^m_{q(n)}$.
\par Then $(q,T)$ is called a homomorphism of $\{ {\cal B}_n, \pi ^n_k, \Lambda \} $
into $\{ {\cal C}_n, \breve{\pi }^n_k, \Upsilon \} $.

\par {\bf 29. Theorem.} {\it There exists a covariant functor
from a category of inverse spectra of microbunles ${\cal SB}$ into a
category of microbundles ${\cal CB}$ induced by the operation $\lim
$.}
\par {\bf Proof.} Let $(q,T)$ be a morphism of an inverse spectrum of microbundles
${\bf S}_1=\{ {\cal B}_n, \pi ^n_k, \Lambda \} $ into ${\bf S}_2=\{
{\cal C}_n, \breve{\pi }^n_k, \Upsilon \} $, where ${\cal B}_n={\cal
B}(A_n,E_n,{\bf F}_n,X_n,i_n,p_n)$ with a left module $X_n$ over a
ring ${\bf F}_n$ for each $n\in \Lambda $; ${\cal C}_n={\cal
B}(C_n,D_n,{\bf G}_n,Y_n,\breve{i}_n,\breve{p}_n)$ with a left
module $Y_n$ over a ring ${\bf G}_n$ for each $n\in \Upsilon $. In
view of Theorem 27 there exist limits ${\cal B}(A,E,{\bf F},X,i,p)=
\lim {\bf S}_1 $ and ${\cal B}(C,D,{\bf G},Y,\breve{i}, \breve{p})=
\lim {\bf S}_2$ of the inverse spectra of microbundles. Put
\par $(29.1)$ $(c_k,d_k,v_k,y_k)=t_k(a_{q(k)},b_{q(k)},u_{q(k)}, x_{q(k)})$ \\
for each $k\in \Upsilon $, where $(a,b,u,x)=((a_k,b_k,u_k,x_k): k\in
\Lambda )$ such that \\ $(\forall k\in \Lambda , ~ a_k\in A_k, ~
b_k\in E_k, ~ u_k\in {\bf F}_k, ~ x_k\in X_k)$.
\par For each $k\le n$ in $\Upsilon $ there exists $m\in \Upsilon $
such that $m\ge k$ and $m\ge n$, since a set $\Upsilon $ is
directed. From $(28.3)$ and $(26.2)$ it follows that
\par $(29.2)$ $(c_k,d_k,v_k,y_k)= \breve{\pi }^n_k(c_n,d_n,v_n,y_n)$.
\par Therefore $(29.1)$ and $(29.2)$ imply that
a limit map $t=\lim (q,T)$ exists from ${\cal B}(A,E,{\bf F},X,i,p)$
into ${\cal B}(C,D,{\bf G},Y,\breve{i}, \breve{p})$ such that \par
$(29.3)$ $t=(t^1,t^2,t^3,t^4)$ with \par $(29.4)$ $t^1: A\to C$, $ ~
t^2: E\to D$, $ ~ t^3: {\bf F}\to {\bf G}$, $ ~ t^4: X\to Y$ \\
where ${\bf G}=\lim \{ {\bf G}_n, \breve{\pi }^{3,n}_{3,k}, \Upsilon
\} $ and $Y=\lim \{ Y_n, \breve{\pi }^{4,n}_{4,k}, \Upsilon \} $.
Then $(28.3)$ and $(26.3)$ imply that
\par $(29.5)$ $t^2\circ i=\breve{i}\circ t^1$ and $\breve{p}\circ
t^2=t^1\circ p$.
\par From the construction of a topology on a limit of an inverse
spectrum (see Subsection 27) it follows that $t$ is continuous.
Therefore, applying $(28.3)$, $(26.4)$ and $(26.5)$ we deduce that
\par $(29.6)$ $\breve{h}\circ t^2\circ h^{-1}(a,ux+wy)=$\par $ t
^3(u)\breve{h}\circ t^2\circ h^{-1}(a,x) + t^3(w) \breve{h}\circ
t^2\circ h^{-1}(a,y)$ and \par
\par $(29.7)$ $\breve{\hat{\pi }}_2\circ \breve{h}\circ t^2=t^4\circ \hat{\pi }_2\circ h$
and  $t^4(ux+wy)=t^3(u)t^4(x)+t^3(w)t^4(y)$ \\
for every $u$ and $w$ in ${\bf F}$, $x$ and $y$ in $X$, $a\in A$,
where $h: V\to U\times X$ is a local homeomorphism for an open
subset $V$ in $E$ corresponding to an open neighborhood $U$ of a
point $a$ in $A$. Therefore $t$ is a continuous homomorphism of
microbundles, since it satisfies the conditions $(29.3)$-$(29.7)$.
Note that $t$ is unique, since $\breve{\pi }_n\circ t = t_n\circ \pi
_{q(n)}$ for each $n\in \Upsilon $, where $\pi _n: {\cal B}(A,E,{\bf
F},X,i,p)\to {\cal B}_n$ with $\pi _n = (\pi _{1,n}, \pi _{2,n}, \pi
_{3,n}, \pi _{4,n})$ for each $n\in \Lambda $ (see Subsection 27).
\par It can be easily verified that a
composition $(q_1\circ q_2, T_2\circ T_1)$ of morphisms $(q_1,T_1):
{\bf S}_1\to {\bf S}_2$ and $(q_2,T_2): {\bf S}_2\to {\bf S}_3$ of
inverse spectra ${\bf S}_j=\{ {\cal B}_{j,k}, ~ \mbox{}^j\pi ^n_k:
k\in \Lambda _j \} $ with $j\in \{ 1, 2, 3 \} $ of microbundles
${\cal B}_{j,k}$ is a morphism from ${\bf S}_1$ into ${\bf S}_3$,
where $T_2\circ T_1=\{ \forall k\in \Lambda _3, ~ t_{2,k}\circ
t_{1,q_2(k)}: {\cal B}_{1,q_1(q_2(k))}\to {\cal B}_{3,k} \} $;
$\Lambda _j$ is a directed set for each $j$. Thus the operation of
taking the limit of an inverse spectrum of microbundles induces a
covariant functor $\lim : {\cal SB}\to {\cal CB}$.

\par {\bf 30. Theorem.} {\it Assume that ${\bf S}=\{ {\cal B}_n, \pi ^n_k, \Lambda \} $
is an inverse spectrum of microbundles and there are homomorphisms
$t_k$ of a microbundle ${\cal C}={\cal B}(C,D,{\bf G},Y,\breve{i},
\breve{p})$ into ${\cal B}_k$ such that $t_k=\pi ^n_k\circ t_n$ for
each $k\le n$ in $\Lambda $. Then there exists a limit homomorphism
$t=\lim \{ t_k: k\in \Lambda \} $, $t: {\cal C}\to {\cal B}(A,E,{\bf
F},X,i,p)$, such that $t_k=\pi _k\circ t$ for each $k\in \Lambda $,
where ${\cal B}(A,E,{\bf F},X,i,p)=\lim \{ {\cal B}_n, \pi ^n_k,
\Lambda \} $. Moreover, if $t_k({\cal C})$ is dense in $X_k$ for
each $k\in \Lambda $, then $t({\cal C})$ is dense in ${\cal
B}(A,E,{\bf F},X,i,p)$.}
\par {\bf Proof.} By virtue of Theorem 29 a continuous homomorphism
$t$ of microbundles exists, since ${\cal C}$ can be written as a
limit of a constant inverse spectrum ${\bf S}_2 = \{ {\cal C}_1, id,
\{ 1 \} \} $ with ${\cal C}_1={\cal C}$ and $\Upsilon = \{ 1 \} $,
where $id$ denotes the identity homomorphism, where the microbundle
${\cal C}$ is given for some left module $Y$ over a topological ring
$\bf G$. Therefore a family $T=\{ t_k: ~ k\in \Lambda \} $ is a
homomorphism from ${\cal C}$ into $\bf S$.
\par Let $(a,e,f,z)\in Q:=\{ (b,d,q,x)\in A\times E\times {\bf F}\times X: ~ p(d)=b \} $.
We take a neighborhood $R:=\{ U\times V\times P\times S: ~ p(V)=U \}
$ of $(a,e,f,z)$, where $U$ is an open subset in a base space $A$, $
~ V$ is an open subset in a total space $E$, $ ~ P$ is an open
subset in a topological ring $\bf F$, $ ~ S$ is an open subset in a
topological left $\bf F$ module $X$. Since $t_k(C,D,{\bf G},Y)$ is
dense in $Q_k:=\{ (a_k,e_k,f_k,z_k)\in A_k\times E_k\times {\bf
F}_k\times X_k: ~ p_k(e_k)=a_k \} $, then an open subset $R_k:=\{
U_k\times V_k \times P_k\times S_k: ~ p_k(V_k)=U_k \} $ in $Q_k$
exists such that $(a,e,f,z)\in \pi _k^{-1} (R_k)$ and $\pi _k^{-1}
(R_k)\subset R$. There exists $(c,d,g,y)\in W := \{ (b,d,q,v)\in
C\times D\times {\bf G}\times Y: ~ \breve{p}(d)=b \} $ such that
$t_k(c,d,g,y)\in R_k$. Hence $t(c,d,g,y)\in \pi _k^{-1} (R_k)$ and
consequently, $t({\cal C})$ is dense in ${\cal B}(A,E,{\bf
F},X,i,p)$.

\par {\bf 31. Proposition.} {\it Let $(id,T):
{\bf S}_1\to {\bf S}_2$ be a homomorphism of inverse spectra ${\bf
S}_j=\{ {\cal B}_{j,k}, ~ \mbox{}^j\pi ^n_k: k\in \Lambda \} $ with
$j\in \{ 1, 2 \} $ of microbundles ${\cal B}_{j,k}$ and let $t_k :
{\cal B}_{1,k}\hookrightarrow {\cal B}_{2,k}$ be an embedding for
each $k\in \Lambda $. Then $\mbox{}^1\pi ^n_k=\mbox{}^2\pi ^n_k\circ
t_k|_{{\cal B}_{1,k}}$ for each $k\le n$ in $\Lambda $ and $\lim
{\bf S}_1=\bigcap_{k\in \Lambda } \mbox{}^2\pi _k^{-1}(t_k({\cal
B}_{1,k}))$.}
\par {\bf Proof.} From the condition $\mbox{}^2\pi ^n_k\circ t_n=t_k\circ
\mbox{ }^1\pi ^n_k$ for each $k\le n\in \Lambda $ and Theorem 29 it
follows that $\mbox{}^1\pi ^n_k=\mbox{}^2\pi ^n_k\circ t_n|_{{\cal
B}_{1,n}}$ for each $k\le n$ in $\Lambda $, since $t_k$ is the
embedding. Then a limit $t:=\lim T$ is a bijective map from ${\cal
B}_1=\lim {\bf S}_1$ into ${\cal B}_2=\lim {\bf S}_2$. For an
element $(a,b,c,x)\in {\cal B}_1$ take a neighborhood $W$. A
neighborhood $V_k$ of $\mbox{}^1\pi _k(a,b,c,x)$ in ${\cal B}_{1,k}$
exists such that $(\mbox{}^1\pi _k)^{-1}(V_k)\subset W$,
consequently, $V:=(\mbox{}^2\pi _k)^{-1}(t_k(V_k))$ is a
neighborhood of $t(a,b,c,x)$ in $t({\cal B}_1)$ such that
$t^{-1}(V)\subset W$. Therefore ${\cal B}_1\subset {\cal B}_2$ and
${\cal B}_1\subset \bigcap_{k\in \Lambda } \mbox{}^2\pi
_k^{-1}(t_k({\cal B}_{1,k}))$. On the other hand, if $(d,e,q,y)\in
\bigcap_{k\in \Lambda } \mbox{}^2\pi _k^{-1}(d_k,e_k,q_k,x_k)$, then
a family $ \{ \mbox{}^2\pi _k(d,e,q,y): ~ k\in \Lambda \} $ is a
threat of the spectrum ${\cal S}_1$, consequently, $\bigcap_{k\in
\Lambda } \mbox{}^2\pi _k^{-1}(t_k({\cal B}_{1,k}))\subset {\cal
B}_1$.

\par {\bf 32. Theorem.} {\it Let $s_1: X_{\bf F}\to X_{\bf G}$ and $s_2: {\bf F}\to {\bf G}$
be homomorphisms of left modules $X_{\bf F}$ and $X_{\bf G}$ and of
rings $\bf F$ and $\bf G$ correspondingly such that \par $(32.1)$
$s_1(v_1x_1+v_2x_2) = s_2(v_1)s_1(x_1) + s_2(v_2)s_1(x_2)$ \\ for
each $v_1$ and $v_2$ in $\bf F$, $ ~ x_1$ and $x_2$ in $X_{\bf F}$.
Let also ${\cal B}(A,E,{\bf F},X_{\bf F},i,p)$ be a microbundle over
$\bf F$. Then a microbundle ${\cal B}(A,E',{\bf G},X_{\bf G},i',p')$
exists such that there is a homomorphism $\pi : {\cal B}(A,E,{\bf
F},X_{\bf F},i,p)\to {\cal B}(A,E',{\bf G},X_{\bf G},i',p')$ with
$s_1\circ \hat{\pi }_2=\hat{\pi }_2'\circ \pi ^2$. Moreover, if
$s_1$ and $s_2$ are either surjective or bijective, then ${\cal
B}(A,E',{\bf G},X_{\bf G},i',p')$ and $\pi$ can be chosen such that
the homomorphism $\pi $ is either surjective or bijective
correspondingly.}
\par {\bf Proof.} Let $e\in E$ and $b\in A$ be such that $i(b)=e$.
We take a neighborhood $V$ of $e$ homeomorphic with $U\times X_{\bf
F}$, where $h$ is a homeomorphism from $V$ onto $U\times X_{\bf F}$,
where $U$ is a neighborhood of $b$ in $A$. We put $\pi _2(V)=V'$ to
be homeomorphic with $U\times X_{\bf G}$ with a homeomorphism $h':
V'\to U\times X_{\bf G}$ and projections $\hat{\pi }_1': V'\to U$, $
~ \hat{\pi }_2': V'\to X_{\bf G}$. This implies that $s_1\circ
\hat{\pi }_2|_V=\hat{\pi }_2'\circ \pi ^2|_V$ (see also the notation
$(29.3)$ and $(29.4)$). Therefore substituting $V$ on $V'$ and $h$
on $h'$ and $X_{\bf F}$ on $X_{\bf G}$ induces maps $i'|_U$,
$p'|_{V'}$, $\iota _0'|_U$ such that $p'|_{V'}\circ i'|_U=\hat{\pi
}_1'\circ \iota _0'|_U$ and $\hat{\pi }_1'\circ h'|_{V'}=p'|_{V'}$.
\par If $U$ and $U_1$ are open neighborhoods of $b$ and $b_1$ in
$A$, $ ~ V=h^{-1}(U\times X_{\bf F})$ and $V_1=h^{-1}(U_1\times
X_{\bf F})$ then $\pi ^2(V)\cap \pi ^2(V_1)=\pi ^2(V\cap
V_1)={h'}^{-1}((U\cap U_1)\times X_{\bf G})$. This provides and
equivalence relation $\Xi $ for each $v\in V'$ and $v_1\in {V_1}'$:
$v\Xi v_1$ if and only if $\hat{\pi }_1'\circ h'(v)=\hat{\pi
}_1'\circ h'(v_1)$ and $\hat{\pi }_2'\circ h'(v)=\hat{\pi }_2'\circ
h'(v_1)$. Using the latter property we choose as a total space
\par $(32.2)$ $E' = \bigcup \{ V': ~ \exists U, ~ U \mbox{ is open in } A, ~ V'=\pi ^2(V), ~
V=h^{-1}(U\times X_{\bf F}) \} / \Xi $. \\
Bases of topologies on $A$ and $X_{\bf G}$ induce a base of a
topology on $E'$. Hence a microbundle ${\cal B}(A,E',{\bf G},X_{\bf
G},i',p')$ exists.
\par Then we put $\pi ^1=id$,
$ ~ \pi ^3=s_2$, $~ \pi ^4 =s_1$ and take a combination \par $\pi
^2= \nabla \{ \pi ^2|_V: ~ \exists U, ~ U \mbox{ is open in } A, ~
V=h^{-1}(U\times X_{\bf F}) \} $. \par From Condition $(32.1)$ it
follows that
\par $(32.3)$ $\pi ^{2}\circ i=i'\circ \pi ^{1}$
and $p'\circ \pi ^{2}=\pi ^{1}\circ p$ and \par $(32.4)$ $h'\circ
\pi ^{2}\circ h^{-1}(b,v_1x_1+v_2x_2)=$\par $ \pi ^{3}(v_1)h'\circ
\pi ^{2}\circ h^{-1}(b,x_1) + \pi ^{3}(v_2) h'\circ \pi ^{2}\circ
h^{-1}(b,x_2)$ and
\par $(32.5)$ ${\hat{\pi }'}_2\circ {h'} \circ \pi ^2 = \pi ^4\circ \hat{\pi }_2\circ h$\\
for every $v_1$ and $v_2$ in ${\bf F}$, $x_1$ and $x_2$ in $X_{\bf
F}$, $b\in A$. Thus from $(32.3)$ and $(32.4)$ this provides a
homomorphism $\pi =(\pi ^1, \pi ^2, \pi ^3, \pi ^4 )$ from ${\cal
B}(A,E,{\bf F},X_{\bf F},i,p)$ into ${\cal B}(A,E',{\bf G},X_{\bf
G},i',p')$.
\par In particular, if $s_1$ and $s_2$ are surjective, $s_1(X_{\bf F})=X_{\bf G}$ and $s_2({\bf F})=\bf G$,
then from Formula $(32.2)$ it follows that $\pi ^2(E)=E'$,
consequently, $\pi $ is surjective. If $s_1$ and $s_2$ are
bijective, then from Identities $(32.3)$-$(32.5)$ we infer that $\pi
^2: E\to E'$ is bijective and hence $\pi $ is bijective.

\par {\bf 33. Corollary.} {\it If Conditions of Theorem 32 are
satisfied and a microbundle ${\cal B}(A,E',{\bf G},X_{\bf G},i',p')$
is provided by this theorem and $s_1$ and $s_2$ are isomorphisms,
then $\pi $ is an isomorphism of a microbundle ${\cal B}(A,E,{\bf
F},X_{\bf F},i,p)$ with ${\cal B}(A,E',{\bf G},X_{\bf G},i',p')$.
Particularly, if $s_1: X_{\bf F}\to X_{\bf F}$ and $s_2: {\bf F} \to
\bf F$ are automorphisms, then $\pi $ is an automorphism of ${\cal
B}(A,E,{\bf F},X_{\bf F},i,p)$.}

\par {\bf 34. Theorem.} {\it Let $A$ be a Tychonoff base space and let
${\cal B}(A,E,{\bf F},X,i,p)$ be a microbundle. Then for each
compactification $cA$ of $A$ there exists a microbundle ${\cal
B}(cA,E',{\bf F},X,i',p')$ and an embedding $\pi : {\cal B}(A,E,{\bf
F},X,i,p)\hookrightarrow {\cal B}(cA,E',{\bf F},X,i',p')$ such that
$i'\circ \pi ^1=\pi ^2\circ i$, $~p'\circ \pi ^2=\pi ^1\circ p$,
$~\pi ^3=I$, $~\pi ^4=I_X$, where $I$ denotes the identity map on a
ring ${\bf F}$, $~ I_X$ is the identity map on a left ${\bf F}$
module $X=X_{\bf F}$, $ ~ \pi = (\pi ^1, \pi ^2, \pi ^3, \pi ^4)$,
$~\pi ^1: A\to cA$, $~\pi ^2: E\to E'$, $E'={E'}_c$. Moreover,
${\cal B}(\beta A,{E'}_{\beta },{\bf F},X,i',p')$ is a maximal
microbundle among such extensions of ${\cal B}(A,E,{\bf F},X,i,p)$,
where $\beta A$ denotes the Stone-$\check{C}$ech compactification of
$A$.}
\par {\bf Proof.} By virtue of Theorem 3.5.1 in \cite{eng}
there exists a compactification $cA$ of $A$. Take an arbitrary fixed
point $b\in A$. There exists a neighborhood $U$ of $b$ in $A$ such
that $i|_U: U\to E$ is continuous with $p\circ i|_U=id|_U$ and
$\iota _0|_U: U\to U\times X_{\bf F}$ is an embedding. Then $h: V\to
U\times X_{\bf F}$ is a homeomorphism, where $V=i(U)$. \par Take
$U'$ open in $cA$ such that $U'\cap c(A)=c(U)$, consequently,
$(U'\times X_{\bf F})\cap (c(A)\times X_{\bf F})=c(U)\times X_{\bf
F}$, where $c: A\hookrightarrow cA$ denotes a homeomorphic
embedding. Therefore, a topological space $V'$ and a homeomorphism
$h': V'\to U'\times X_{\bf F}$ exist such that ${h'}^{-1}\circ
(c\times I_X)|_{U\times X_{\bf F}}=h^{-1}|_{U\times X_{\bf F}}$,
where $I_X$ denotes an identity map on $X_{\bf F}$. Evidently, there
are extensions $i': U'\to V'$ and $p': V'\to U'$ and $\iota _0':
U'\to U'\times X_{\bf F}$ and $\hat{\pi }_1': U\times X_{\bf F}\to
U'$ and $\hat{\pi }_2': U'\times X_{\bf F}\to X_{\bf F}$ such that
$p'\circ i'|_{U'}=id|_{U'}$ and $\hat{\pi }_1'\circ \iota
_0'|_{U'}=p'|_{V'}\circ i'|_{U'}$, where $\iota _0'(b)=(b,0)$ for
each $b\in U'$, $\hat{\pi }_1'(b,x)=b$ and $\hat{\pi }_2'(b,x)=x$
for each $(b,x)\in U'\times X_{\bf F}$.
\par Then $V'\cap V_1'={h'}^{-1}|_{(U'\cap U_1')\times X_{\bf F}}$,
where $U_1$ is an open neighborhood of $b_1\in A$ such that
$i|_{U_1}: U_1\to E$ is continuous with $p\circ i|_{U_1}=id|_{U_1}$
and $\iota _0|_{U_1}: U_1\to U_1\times X_{\bf F}$ is an embedding.
This induces an equivalence relation $v\Xi _1v_1$ for each $v\in V'$
and $v_1\in V_1'$ if and only if $\hat{\pi }_1'\circ h'(v)=\hat{\pi
}_1'\circ h'(v_1)$ and $\hat{\pi }_2'\circ h'(v)=\hat{\pi }_2'\circ
h'(v_1)$.  We put
\par $(34.1)$ $E' = \bigcup \{ V': ~ \exists U', ~ U' \mbox{ is open
in } cA, ~ V'={h'}^{-1}(U'\times X_{\bf F}) \} / \Xi _1$.
\par Therefore there exists an embedding $\pi ^2: E\to E'$ such that
$h'\circ \pi ^2(V) = c(U)\times X_{\bf F}$ for each open $U$ and $V$
as described above. Taking $\pi ^1=c$ and $\pi ^3=I$ and $\pi
^4=I_X$, we deduce from $(34.1)$ and the construction above that
$i'\circ \pi ^1=\pi ^2\circ i$, $~p'\circ \pi ^2=\pi ^1\circ p$.
\par In view of Theorem 3.6.1 in \cite{eng} there is an embedding
of \par ${\cal B}(cA,{E'}_c,{\bf F},X,i',p')\hookrightarrow {\cal
B}(\beta A,{E'}_{\beta },{\bf F},X,i',p')$ \\ for each
compactification $cA$ of $A$, consequently, ${\cal B}(\beta
A,{E'}_{\beta },{\bf F},X,i',p')$ is maximal among such extensions
of the microbundle ${\cal B}(A,E,{\bf F},X,i,p)$.

\par {\bf 35. Corollary.} {\it Suppose that the conditions of Theorem 34
are satisfied and a ring $\bf F$ and a left module $X_{\bf F}=X$
have compactifications $c_{3,1}\bf F$ and $c_{3,2}X_{\bf F}$ which
are a ring $\bf G$ and a left module isomorphic with a left module
$X_{\bf G}$ over $\bf G$ respectively. Then there exists a
microbundle ${\cal B}(A',E',{\bf G},X_{\bf G},i',p')$ such that
$A'=c_1A$ and $E'=c_2E$ are compactifications of $A$ and $E$
respectively and an embedding ${\bf c}: {\cal B}(A,E,{\bf F},X_{\bf
F},i,p)\hookrightarrow {\cal B}(A',E',{\bf G},X_{\bf G},i',p')$
exists.}
\par {\bf Proof.} A procedure of taking an extension ${\cal B}_1$ of the microbundle
\\ ${\cal B}(A,E,{\bf F},X_{\bf F},i,p)$ at first by Theorem 32 and
then extending ${\cal B}_1$ by Theorem 34 provides the microbundle
${\cal B}(A',E',{\bf G},X_{\bf G},i',p')$ with a compact base space
$A'$. From Formulas $(32.2)$ and $(34.1)$ we deduce that a total
space $E'$ is compact, since $A'$, $ ~ \bf G$ and $X_{\bf G}$ are
compact.

\par {\bf 37. Definition.} Let ${\cal S} = \{ {\cal B}_j, \pi ^k_j, J \} $ be
an inverse spectrum of microbundles ${\cal B}_j={\cal
B}(A_j,E_j,i_j,p_j)$, and let \par $\lim \{ \pi ^k_j: ~ j\in J, ~
k\in J, ~ j\le k < l \} : ~ {\cal B}_l\to \lim \{ {\cal B}_j, \pi
^k_j, J: j\le k<l \} $ \\ be a homeomorphism for each limit element
$l\in J$, then ${\cal S}$ is called continuous.
\par Let $\tau $ be an infinite cardinal. A directed set $J$ is
called $\tau $-complete, if each its linearly ordered subset $K$ of
the cardinality $card (K)\le \tau $ has a supremum in $J$.
\par If an inverse spectrum ${\cal S}$ of microbundles ${\cal B}_j={\cal
B}(A_j,E_j,{\bf F}_j,X_j,i_j,p_j)$ with a $\tau $-complete directed
set $J$ is continuous and there is a least element $j_0$ in $J$ and
$\pi ^k_j({\cal B}_k)={\cal B}_j$ with a compact base space $A_j$
for each $j\le k\in J$, then $\cal S$ is called $\tau $-complete.
\par Assume that an inverse spectrum ${\cal S}$ of microbundles
is $\tau $-complete and $wA_j\le \tau $ for each $j\in J$, where
$wA_j$ denotes a weight of $A_j$, then $\cal S$ is called a $\tau
$-spectrum.

\par {\bf 38. Proposition.} {\it Let ${\cal S}_1=
\{ {\cal B}_{1,k}, ~ \mbox{}^1\pi ^n_k: k\in \Lambda \} $ be a
continuous $\tau $-complete spectrum of microbundles and ${\cal
S}_2=\{ {\cal B}_{2,k}, ~ \mbox{}^2\pi ^n_k: k\in \Lambda \} $ be a
$\tau $-spectrum of microbundles, ${\cal B}_j=\lim {\cal S}_j$,
$f_{j_0}: {\cal B}_{1,j_0}\to {\cal B}_{2,j_0}$ and $f: {\cal
B}_1\to {\cal B}_2$ be homomorphisms of microbundles such that $\pi
_{2,j_0}\circ f= f_{j_0}\circ \pi _{1,j_0}$. Let also either
$E_{j,k}$ be compact for each $j$ and $k$, or ${\bf F}_{1,k}={\bf
F}_{2,k}$ and $X_{1,k}=X_{2,k}$ for each $k$. Then $f$ is a limit of
homomorphisms between cofinal subspectra of ${\cal S}_1$ and ${\cal
S}_2$.}
\par {\bf Proof.} According to the conditions of this proposition
$\mbox{}^j\pi ^{1,n}_{1,k}(A_{j,n})=A_{j,k}$ and $\mbox{}^j\pi
^{2,n}_{2,k}(E_{j,n})=E_{j,k}$ and $\mbox{}^j\pi ^{3,n}_{3,k}({\bf
F}_{j,n})={\bf F}_{j,k}$ and $\mbox{}^j\pi
^{4,n}_{4,k}(X_{j,n})=X_{j,k}$ for each $k\le n$ in $J$ and $j\in \{
1, 2 \} $. Then $f=(f^1,f^2,f^3,f^4)$, where $f^1: A_1\to A_2$,
$f^2: E_1\to E_2$, $f^3: {\bf F}_1\to {\bf F}_2$, $f^4: X_1\to X_2$,
where a microbundle ${\cal B}_j={\cal B}(A_j,E_j,{\bf
F}_j,X_j,i_j,p_j)$ is on a left module $X_{{\bf F}_j}=X_j$ over a
ring ${\bf F}_j$ for each $j \in \{ 1, 2 \} $. If $E_{j,k}$ is a
compact total space, then a ring ${\bf F}_{j,k}$ and a left module
$X_{j,k}$ are compact as follows from Theorem 3.1.10 in \cite{eng}
and Condition $(2.3)$ above. By virtue of Theorem 40 in
\cite{scepumn81} (see also \cite{fedumn86}) and Theorems 29 and 34
and Proposition 31 above there exists a cofinal subset $\Lambda $ in
$J$ such that $f=\lim \{ t_k: ~ k\in \Lambda \} $, where $t_k: {\cal
B}_{1,k}\to {\cal B}_{2,k}$ is a homomorphism for each $k\in \Lambda
$.
\par {\bf 39. Corollary.} {\it If the conditions of Theorem 38
are satisfied and $f_0$ and $f$ are homeomorphisms, then $f$ is a
limit of homeomorphisms between cofinal subspectra of ${\cal S}_1$
and ${\cal S}_2$.}

\end{document}